\documentclass{article}
 \usepackage{amsmath}

\newtheorem{thm}{Theorem}[section]
\newtheorem{cor}[thm]{Corollary}
\newtheorem{prop}[thm]{Proposition}

\newtheorem{lem}[thm]{Lemma}
\newtheorem{Def}[thm]{Definition}
\newtheorem{rem}[thm]{Remark}

\newtheorem{ex}[thm]{Example}

\newcommand{\be}{\begin{equation}}
\newcommand{\ee}{\end{equation}}
\newcommand{\ben}{\begin{enumerate}}
\newcommand{\een}{\end{enumerate}}
\newcommand{\beq}{\begin{eqnarray}}
\newcommand{\eeq}{\end{eqnarray}}
\newcommand{\beqn}{\begin{eqnarray*}}
\newcommand{\eeqn}{\end{eqnarray*}}

\newcommand{\pa}{\partial}

\newcommand{\qed}{\hspace*{\fill}Q.E.D.}  

\begin{document}
\title{On Geodesics of Sprays and Projective Completeness}
\author{Guojun Yang}
\date{}
\maketitle
\begin{abstract}
 Geodesics, which  play an important role in spray-Finsler geometry, are integral curves of
  a spray vector field on a manifold.
 Some comparison theorems and rigidity issues are established on the
 completeness of geodesics of a spray or a Finsler metric. In this
 paper, projectively flat sprays with weak Ricci constant (eps. constant
 curvature)  are classified at the level of geodesics. Further, a geodesic method is
 introduced to determine an $n$-dimensional spray based on a family of
 curves with $2(n-1)$ free constant parameters as geodesics.
 Finally, it shows that a spray is projectively complete under certain
 condition satisfied by the domain of geodesic parameter of all
 geodesics.


 {\bf MR(2000) subject classification: }
53B40, 53C60
\end{abstract}

\section{Introduction}

Spray geometry  studies the properties of  sprays on a manifold,
and it is closely related to Finsler geometry. Every Finsler
metric induces a natural spray but there are a lot of sprays which
are not Finsler-metrizable (not be induced by any Finsler metric)
(\cite{EM, LMY, Yang1}). So a popular topic is to investigate
whether a given spray is metrizable or not, and what's more
important is to give necessary and sufficient conditions for
certain class of sprays to be metrizable (\cite{BM2, Yang2,
Yang3}). It is also important to investigate the properties of
some special classes of sprays, for example, (locally)
projectively flat sprays, Berwald sprays, sprays of scalar (resp.
isotropic, constant) curvature, Hamel (resp. Funk) sprays
(\cite{BM2, LS, Yang2, Yang3}).

A spray ${\bf G}$ on a manifold $M$ defines a special vector filed
on a conical region $\mathcal{C}$ of
 $TM\setminus \{0\}$, and  it naturally defines its integral curves and the
 projections of the integral curves onto the manifold $M$ are
 called geodesics. Geodesics play an important role in the studies
 of comparison theorems and rigidity issues  on  spray or Finsler
 manifolds. In \cite{Shen3},  Z. Shen studies two
 pointwise projectively related Einstein Finsler metrics and determine the metrics along geodesics.
  In \cite{Yang1}, the present author obtains a comparison theorem on the
 Ricci curvatures of a spay and a Finsler metric which are
 pointwise projectively related and the corresponding projective
 factor is estimated. In \cite{Br}, R. Bryant proves that a
 geodesically reversible Finlser metric on $S^2$ with positive
 constant flag curvature is a Riemann metric. In \cite{Ro}, C.
 Robles classifies geodesics of Randers metrics of constant
flag curvature. In \cite{HM}, L. Huang and X. Mo obtain the
relation between the geodesics of two Finsler metrics $F$ and
$\tilde{F}$, where $\tilde{F}$ is defined by the navigation data
$(F,V)$ with $V$ being a homothetic vector field of $F$. In this
paper, we study projectively flat sprays with weak Ricci constant,
 the construction of
 sprays from a geodesic method and  the projective completeness of sprays.

In \cite{Yang2}, it introduces  sprays of constant curvature and a
spray {\bf G} of constant curvature is weakly Ricci constant (the
Ricci curvature is constant along any geodesic of {\bf G}). For
two pointwise projectively related sprays, they have same
geodesics as point sets and their geodesic parameters are closely
related by the projective factor. Starting from this fact, we can
determine a projectively flat spray with  weak Ricci constant at
geodesic level.

 We consider a projectively flat spray manifold
$({\bf G}, M)$, that is,
 \be\label{Geody2}
 G^i=\widetilde{G}^i+Py^i,
 \ee
where ${\bf \widetilde{G}}$ is a locally Minkowski spray on $M$.
We have the following theorem.

\begin{thm}\label{PGeod1}
  If the spray  {\bf G} in (\ref{Geody2}) is weakly Ricci constant $Ric_{;0}=0$ or of constant curvature, then
 along any geodesic $x=x(s)$ of {\bf G}, $P(s):=P(x(s),x'(s))$ is
 given by one of the following cases:
  \be\label{Geod003}
 P(s)=\frac{1}{s+\kappa},\ \ P(s)=-c\cdot \tan(cs+\kappa),\ \
 P(s)=-\frac{c(1-\kappa e^{2cs})}{1+\kappa e^{2cs}},
  \ee
  where $c,\kappa$ are constant. Further, if {\bf G} is complete,
  then $P(s)$ is given by
   \be\label{Geod004}
P(s)=-\frac{c(1-\kappa e^{2cs})}{1+\kappa e^{2cs}}.
   \ee
\end{thm}

In Theorem \ref{PGeod1}, we can further give the relation between
the geodesic parameters of {\bf G} and $\widetilde{G}$ by
(\ref{Geod003}) (see Proposition \ref{PGeod2}, Corollary
\ref{CGeod1}).

The family of geodesics of an $n$-dimensional  spray considered as
point sets or paths is dependent on $2(n-1)$ free constant
parameters. A path space is a family of curves satisfying certain
conditions (Definition \ref{DGeod00}). We can freely give many
interesting path spaces, especially in dimension two. Starting
from a path space, we can construct its corresponding spray.

\begin{thm}\label{LGeody3}
 In an  $n$-dimensional path space $\mathcal{G}$, all paths in a local coordinate system $(x^i)$ can be  parameterized
 under  a  variable $t$
with $2(n-1)$ free constant parameters $u,v$ as follows:
 \be\label{Geodj2}
 x=x(t)=\sigma(t;u,v),\ \ \ \ (u,v\in R^{n-1}).
 \ee
 Further, the parametric equation (\ref{Geodj2})
 induces a spray {\bf G} whose  geodesics  are
 given by (\ref{Geodj2}) with
 $t$ as its geodesic parameter, and if  a new
 variable $s=s(t)=s(t;u,v)$ is given with $s'(t)>0$, then  it gives a spray ${\bf \bar{G}}\in Proj({\bf {G}})$
 with $s$ as its geodesic parameter.
 \end{thm}

If a family of curves can be parameterized in the form
(\ref{Geodj2}), then with an auxiliary parameter $c>0$ multiplied
 by $t$ in (\ref{Geodj2}), we can
obtain the corresponding spray by eliminating the parameters
$u,v,c,t$. We give some examples to show how to solve the sprays
from given path spaces (see Examples
\ref{EGeody001}-\ref{EGeody004}).

In  the study of rigidity issues on a Finsler or spray manifold,
it is important to assume that the (Finsler) spray in
consideration be (positively/negatively) complete. A given spray
is not necessarily (positively/negatively) complete. So a natural
problem is whether a spray can be projectively
(positively/negatively) complete or not. We solve this problem
under certain conditions in the following result (Theorem
\ref{th3}).

\begin{thm}\label{th3}
 Let {\bf G} be a spray on a manifold $M$ with its each geodesic $x=x(t)$ being defined on the maximal interval
 $I$ given by one case of the following
 \be\label{Geodc1}
 I= (a,b),\ \ or\ \ (a,+\infty), \ \ or\ \
(-\infty,b),
 \ee
where $a=a(u,v)<0, b=b(u,v)>0$ with $u=x(0), v=x'(0)$ are
$C^{\infty}$
 functions on a conical region $\mathcal{C}$ of $TM\setminus \{0\}$.
  Then {\bf G} is projectively (positively/negatively)
 complete on $\mathcal{C}$.
\end{thm}

 In Theorem \ref{th3}, usually we can also put $u,v$ as that in (\ref{Geodj2}) (see
 Example \ref{EGeod4}). If (\ref{Geodc1}) is not satisfied, it is
 uncertain that {\bf G} is projectively complete (cf. Example \ref{EGeod4})).
We give  Examples \ref{EGeod1}-\ref{EGeod4} as an application of
Theorem \ref{th3}. A Finsler metric is not necessarily
projectively (positively/negatively) complete, namely, if {\bf G}
in Theorem \ref{th3} is a Finsler spray,  the
 corresponding spray projective to {\bf G} may not be a Finsler spray.

\section{Geodesic parameters in projective relations}

 A {\it spray} on $M$, in our consideration, is a smooth
vector field ${\bf G}$ on a conical region $\mathcal{C}$ of
 $TM\setminus \{0\}$ (an important case is $\mathcal{C}=TM\setminus \{0\}$) expressed in
a local coordinate system $(x^i,y^i)$ in $TM$ as follows
 $${\bf G}=y^i\frac{\pa}{\pa x^i}-2G^i\frac{\pa}{\pa y^i},$$
 where $G^i$ are local homogeneous functions satisfying
 $G^i(x,\lambda y)=\lambda^2G^i(x,y)$ for $\lambda>0$. If
 $\mathcal{C}=TM\setminus \{0\}$, {\bf G} is called {\it regular};
 otherwise, it is called {\it singular}.

The integral curves of {\bf G} projected onto $M$ are the
geodesics of {\bf G}. Let $x=x(s)$ be a geodesic of {\bf G}. Then
it satisfies the following ODE:
 $$
\frac{d^2x^i}{ds^2}+2G^i(x,\frac{dx}{ds})=0,
 $$
where $s$ is called a {\it geodesic parameter} of the geodesic
$x=x(s)$. Reparameterizing a geodesic $x=x(s)$ by a general
parameter $t$ with $ds/dt>0$, we have
 \be\label{S2}
\frac{d^2x^i}{dt^2}+2G^i(x,\frac{dx}{dt})=\gamma(t)\frac{dx^i}{dt},
 \ee
where $\gamma(t)$ is given by
 \be\label{S3}
\gamma(t)=\frac{d^2s}{dt^2}\Big/\frac{ds}{dt}=-\frac{d^2t}{ds^2}\Big/(\frac{dt}{ds})^2.
 \ee

Let ${\bf G}, \bar{{\bf G}}$ be two sprays pointwise projectively
related by $\bar{G}^i=G^i+Py^i$. Let $x=x(t)$ be a geodesic of
{\bf G} or $\bar{{\bf G}}$ as a point set for a general parameter
$t$. Then along the geodesic $x=x(t)$, it follows from (\ref{S2})
and (\ref{S3}) that
 \be\label{Geod1}
2P(t)= \frac{\bar{s}''(t)}{\bar{s}'(t)}-\frac{s''(t)}{s'(t)},\ \ \
\ \big(P(t):=P(x(t),x'(t))\big),
 \ee
where $s,\bar{s}$ are the geodesic parameters of the curve
$x=x(t)$ in ${\bf G}, \bar{{\bf G}}$ respectively. In particular,
along a geodesic $x=x(s)$ of {\bf G}, it follows from
(\ref{Geod1}) that
 \be\label{Geod2}
2P(s)= \frac{\bar{s}''(s)}{\bar{s}'(s)}, \ \ \ \
\big(P(s):=P(x(s),x'(s))\big),
 \ee
If we express the geodesic $x=x(s)$ of {\bf G} as the geodesic
$x=x(\bar{s})$ of ${\bf \bar{G}}$,  by (\ref{Geod2}), we have
 \be\label{Geody1}
 2P(\bar{s})=2P(x(s),x'(s))\frac{ds}{d\bar{s}}=\frac{\bar{s}''(s)}{\big(\bar{s}'(s)\big)^2},
 \ \ \ \ \big(P(\bar{s}):=P(x(\bar{s}),x'(\bar{s}))\big).
 \ee
So if  $P(s)$ or $P(\bar{s})$ is known, the relation
$\bar{s}=\bar{s}(s)$ can be obtained  from (\ref{Geod2}) or
(\ref{Geody1}).

\begin{ex}\label{EGeod01}
 Let $F$ be the Funk metric on a strongly convex domain $\Omega\subset R^n$. Define a projectively flat spray {\bf G} by
 $$
G^i=Py^i,\ \ \ P:=cF,
 $$
 where $c$ is a constant. Any geodesic $x=x(t)$ (as a point set)
 of {\bf G} is given by
  $$x=x(t)=vt+u,\ \ \  \big(-\frac{1}{F(u,-v)}<t<\frac{1}{F(u,v)}\big),$$
 where $u,v\in R^n$ are constant vectors. We have
  \be\label{cwy3}
 F(vt+u,v)=\frac{F(u,v)}{1-tF(u,v)}.
  \ee
 Let $s$ be a geodesic
 parameter of {\bf G}. Then by (\ref{Geod2}) and (\ref{cwy3}) we have
  \be\label{Geodg1}
 \frac{s''(t)}{s'(t)}=2cF(vt+u,v)=\frac{2cF(u,v)}{1-tF(u,v)},
  \ee
  integration of which with $s(0)=0$ gives
 \be\label{Geod002}
   s= s(t)=
 \begin{cases}
  \kappa\ln\big[1-tF(u,v)\big],\hspace{2.15cm}  (c=\frac{1}{2}),\\
\kappa \Big[1-\big(1-tF(u,v)\big)^{1-2c}\Big], \hspace{1cm} (c\ne
\frac{1}{2}),
 \end{cases}
 \ee
 where $\kappa$ is a constant with $\kappa< 0$ for  $c\ge 1/2$,
 and $\kappa> 0$ for  $c< 1/2$.
     Thus the spray is positively complete for  $c\ge 1/2$, and any geodesic is defined on a finite open interval for $c<1/2$.
     Besides, the spray {\bf G} is (locally) metrizable if and
     only if $c=0,1,1/2$ (see \cite{Yang1}).
\end{ex}

\begin{ex}\label{EGeod001}
 In Example \ref{EGeod01}, if the spray {\bf G} is given by
  $$
 G^i(y):=Py^i,\ \ \ P:=c\big[F(y)-F(-y)\big],
  $$
  then by (\ref{cwy3}) and
  \be\label{cwy003}
F(vt+u,-v)=\frac{F(u,-v)}{1+tF(u,-v)}.
  \ee
  it follows from (\ref{Geod2}) that
   \beqn
\frac{s''(t)}{s'(t)}=\frac{2cF(u,v)}{1-tF(u,v)}-\frac{2cF(u,-v)}{1+tF(u,-v)},
   \eeqn
   integration of which  with $s(0)=0$ gives
   \be\label{Geodg2}
 s=s(t)=\kappa\int^t_0\Big[\big(1-tF(u,v)\big)\big(1+tF(u,-v)\big)\Big]^{-2c}dt,
   \ee
   where $\kappa>0$ is constant. From (\ref{Geodg2}), it is clear to
   conclude that {\bf G} is complete if $c\ge 1/2$; $s$ is bounded
   in a finite open interval if $c<1/2$.
\end{ex}

\section{Projective flat sprays with weak Ricci constant}

For a spray {\bf G}, the  Riemann curvature tensor $R^i_{\ k}$  is
defined by
 $$
  R^i_{\ k}:=2\pa_k G^i-y^j(\pa_j\dot{\pa_k}G^i)+2G^j(\dot{\pa_j}\dot{\pa_k}G^i)-(\dot{\pa_j}G^i)(\dot{\pa_k}G^j),
 $$
where we define $\pa_k:=\pa/\pa x^k,\dot{\pa}_k:=\pa/\pa y^k$. The
trace of $R^i_{\ k}$ is called the Ricci curvature, $ Ric:=R^i_{\
i}$. For a spray tensor $T=T_idx^i$ as an example,   the
horizontal and vertical derivatives of $T$ with respect to Berwald
 connection are given by
  $$
 T_{i;j}=\delta_jT_i-T_rG^r_{ij},\ \ \ \ \ \ \
 T_{i.j}=\dot{\pa}_jT_i,\ \ \ \ (\delta_i:=\pa_i-G^r_i\dot{\pa}_r,\ G^k_{ir}:=\dot{\pa}_r\dot{\pa}_iG^k)).
  $$
A spay is called weakly Ricci constant if
$Ric_{;0}:=Ric_{;r}y^r=0$. A spray ${\bf G}$ is said to be of
 {\it constant curvature} if $R^i_{\ k}$ is given by $R^i_{\ k}=R\delta^i_k-\tau_ky^i$  with
 (\cite{Yang2})
  $$
 \tau_{i;k}=0 \ (\ \Leftrightarrow \ R=\tau_k=0,\ or \ R_{;i}=0 (R\ne 0)).
  $$
By definition, it is clear that a spray of constant curvature is
weakly Ricci constant.   For two pointwise projectively related
sprays ${\bf G}, {\bf\bar{G}}$ with $\bar{G}^i=G^i+Py^i$, their
Ricci curvatures $Ric,\bar{R}ic$ are related by
 \be\label{pr13}
 \bar{R}ic=Ric-(n-1)(P_{;0}-P^2).
 \ee

\

We consider a projectively flat spray manifold $({\bf G}, M)$
given by (\ref{Geody2}), that is,
 $$
 G^i=\widetilde{G}^i+Py^i,
 $$
where ${\bf \widetilde{G}}$ is a locally Minkowski spray on $M$
(${\bf \widetilde{G}}$ has local straight lines as geodesics). If
{\bf G} is weakly Ricci constant, then we can determine the
projective factor $P$ along geodesics, which is shown in Theorem
\ref{PGeod1}.

\

{\it Proof of Theorem \ref{PGeod1} :} By (\ref{pr13}) and
$\widetilde{G}^i=G^i-Py^i$, the Ricci curvature $Ric$ of {\bf G}
is given by
 $$
 Ric=-(n-1)(P^2+P_{;0}).
 $$
Therefore, $Ric_{;0}=0$ is equivalent to $P_{;0;0}+2PP_{;0}=0$.
Then along a geodesic $x=x(s)$ of {\bf G}, we have
 $$
 P''(s)+2P(s)P'(s)=0.
 $$
whose solution is given by one of the three cases in
(\ref{Geod003}). Further, if {\bf G} is complete, it is clear that
(\ref{Geod004}) follows from (\ref{Geod003}).   \qed

\

If the spray  {\bf G} in (\ref{Geody2}) is weakly Ricci constant
$Ric_{;0}=0$, then applying (\ref{Geod003}) and (\ref{Geody1}), we
obtain the following proposition.

\begin{prop}\label{PGeod2}
  Let the spray  {\bf G} in (\ref{Geody2}) be weakly Ricci constant (esp. of
constant curvature).  For any geodesic $\sigma$, let $s$ and $t$
be the geodesic parameters of $\sigma$ with respect to {\bf G} and
${\bf \widetilde{G}}$ respectively. Then $s=s(t)$ is given by one
of the following cases:
 \beq
 &&s=at, (a>0);\hspace{2.43cm} s=b\ln(1+at), (ab>0);
\label{Geody3}\\
 && s=\frac{bt}{1+at}, (a\ne0,b>0);\hspace{0.8cm} s=c\big[\arctan(at+b)-\arctan b\big],(ac>0);\label{Geody4}\\
 &&s=c\ln\frac{1+bt}{1+at}, \big((b-a)c>0,\  ab\ne0\big),\label{Geody04}
 \eeq
 where $a,b,c$ are constant, and in (\ref{Geody04}), it further requires
 $s'(t)>0$ (see Remark \ref{RGeod1}).
\end{prop}

{\it Proof :}  By (\ref{Geody1}) we need to solve the following
ODE with initial conditions:
 $$
 \frac{s''(t)}{s'(t)}=2P(s)s'(t),\ \ \ \ (s(0)=0,\ s'(t)>0),
 $$
integration of which gives
 \be\label{Geody5}
 s'(t)=a e^{2\int P(s)ds},\ \ \ \int  e^{-2\int P(s)ds}ds=at+b,
 \ee
where $a,b$ are two constants. Now $P(s)$ is given by
(\ref{Geod003}) from Theorem \ref{PGeod1}, and thus we can obtain
$s=s(t)$ by plugging $P(s)$ into (\ref{Geody5}).

If $P(s)=0$, then (\ref{Geody5}) gives $s=at+b$. Since
$s(0)=0,s'(t)>0$, we obtain $s=at\ (a>0)$, which gives the first
formula in (\ref{Geody3}).

If $P(s)=c\ne0$ is constant, then (\ref{Geody5}) gives
  the second formula in (\ref{Geody3}) with
$ab>0$.

If $P(s)$ is given by the first formula in (\ref{Geod003}), then
(\ref{Geody5}) gives
 $$
 s=-\kappa+\frac{1}{at+b},
 $$
which can be rewritten as the form of the first formula in
(\ref{Geody4}) by $s(0)=0,s'(t)>0$.

If $P(s)$ is given by the second formula in (\ref{Geod003}) ($c\ne
0$), then (\ref{Geody5}) gives
 $$
 s=-\frac{\kappa-\arctan(at+b)}{c},
 $$
which can be rewritten as the  second formula in (\ref{Geody4}) by
$s(0)=0,s'(t)>0$.

If $P(s)$ is given by the third formula in (\ref{Geod003})
($c\kappa\ne 0$), then (\ref{Geody5}) gives
 $$
s=\frac{1}{2c}\ln\big(\frac{1}{at+b}-\frac{1}{\kappa}\big)
 $$
which can be rewritten as the   formula in (\ref{Geody04}) by
$s(0)=0,s'(t)>0$.     \qed

\

In (\ref{Geody04}), by $s'(t)>0$, we have further restriction on
the constant parameters $a,b,c$, which is shown  in the following
remark.

\begin{rem}\label{RGeod1}
 In (\ref{Geody04}), let $t$ be defined on the maximal interval
 $(\kappa_1,\kappa_2)$ with $\kappa_1<0<\kappa_2$. It is easy to conclude the following
 cases from $s'(t)>0$:
  \beqn
 &&a>0,\ b>0:\
  \begin{cases}
  t\in(\kappa_1,\kappa_2)\subset (-\frac{1}{a},+\infty), \hspace{0.5cm} (b<a)\\
  t\in(\kappa_1,\kappa_2)\subset (-\frac{1}{b},+\infty), \hspace{0.2cm}
  (b>a),
  \end{cases}\\
  &&a<0,\ b<0:\
  \begin{cases}
  t\in(\kappa_1,\kappa_2)\subset (-\infty,-\frac{1}{a}), \hspace{0.5cm} (b>a)\\
  t\in(\kappa_1,\kappa_2)\subset (-\infty,-\frac{1}{b}), \hspace{0.2cm}
  (b<a),
  \end{cases}\\
  &&a>0,\ b<0:\
  t\in(\kappa_1,\kappa_2)\subset(-\frac{1}{a},-\frac{1}{b}), \\
  &&a<0,\ b>0:\
  t\in(\kappa_1,\kappa_2)\subset  (-\frac{1}{b},-\frac{1}{a}).
  \eeqn
\end{rem}

By Proposition \ref{PGeod2} and Remark \ref{RGeod1}, we directly
obtain the following corollary.

\begin{cor}\label{CGeod1}
If the spray {\bf G} in Proposition \ref{PGeod2} ($P\ne 0$) is
complete, then $s=s(t)$ is given by one of the following two
cases:
 \beq
 &&s=b\ln(1+at),\ (ab>0),\label{Geody6}\\
&&s=c\ln\frac{1+bt}{1+at},\ \big((b-a)c>0,\
ab<0\big),\label{Geody7}
 \eeq
 where in (\ref{Geody6}) and (\ref{Geody7}), we respectively have
  \beqn
 &&\hspace{1cm} t\in (-\infty,-\frac{1}{a}) \ \ {\rm if}\ a<0, \ {\rm and} \ t\in
 (-\frac{1}{a},+\infty)\ \ {\rm if}\ a>0;\\
 &&t\in(-\frac{1}{a},-\frac{1}{b})\ \ {\rm if} \ a>0, b<0,\ {\rm and} \
 t\in(-\frac{1}{b},-\frac{1}{a})\ \ {\rm if} \ a<0, b>0.
  \eeqn
\end{cor}

Now in the following, we give some projectively flat sprays to
verify the above results on the projective factors and the
geodesic parameters.

\begin{ex}\label{EGeod02}
Consider the spray {\bf G} in Example \ref{EGeod01}. A direct
computation shows that {\bf G} is weakly Ricci constant or of
constant curvature if and only if $c=0,1,1/2$.  Let $x=x(t)=vt+u$
be a geodesic (as a point set) of {\bf G}, and the geodesic
parameter $s$ in {\bf G} is given by (\ref{Geod002}). Then it
follows from (\ref{Geody1}) and (\ref{Geod002}) that $P=cF$ is
given by
 \be
 P(s)=
  \begin{cases}
  -\frac{1}{2\kappa},\hspace{1.85cm}  (c=\frac{1}{2})\\
  \frac{c}{2c-1}\cdot \frac{1}{s-\kappa},\hspace{1cm}
  (c\ne\frac{1}{2}).
  \end{cases}\label{Geod005}
 \ee
 It is clear from
(\ref{Geod005}) that $P(s)$ is in one of  the forms in
(\ref{Geod003}) if and only if $c=0,1,1/2$. Meanwhile,  $s=s(t)$
is given in Proposition \ref{PGeod2} if and only if $c=1,1/2,1$,
and in this case, $s=s(t)$ is in the respective forms shown in
(\ref{Geody3}) and the first formula in (\ref{Geody4}).
\end{ex}

\begin{ex}
 Let {\bf G} be the spray in Example \ref{EGeod001} with $c=1/2$. {\bf G} is complete and it is of constant curvature.
 Then it follows from (\ref{Geodg2}) that
  \be\label{Geody8}
 s=\kappa\ln\frac{1+tF(u,-v)}{1-tF(u,v)},
  \ee
  where $\kappa>0$ is a constant. In this case, $s=s(t)$ is in the
  form (\ref{Geody7}), and it is easy to verify that $P(s)$ is in the form of the third
  formula in (\ref{Geod003}) by plugging (\ref{Geody8}) into (\ref{Geody1}).
\end{ex}

\begin{ex}
 Let {\bf G} be a spray on $R^n$ defined by
  $$
 G^i:=Py^i,\ \ \ P:=-\frac{\langle x,y\rangle}{1+|x|^2}.
  $$
 {\bf G} is metrizable and it is of constant curvature. Let $x=x(t)=vt+u$
be a geodesic (as a point set) of {\bf G}, and by (\ref{Geod2}),
the geodesic parameter $s$ of {\bf G} satisfies
 $$
 \frac{s''(t)}{s'(t)}=-\frac{2(|v|^2t+\langle
 u,v\rangle)}{1+|u|^2+2\langle u,v\rangle t+|v|^2t^2}.
 $$
 Solving the ODE, we obtain
  $$
 s=s(t)=\kappa_1+\kappa_2 \arctan\frac{|v|^2t+\langle u,v\rangle}{\sqrt{(1+|u|^2)|v|^2-\langle
 u,v\rangle^2}},
  $$
  where $\kappa_1,\kappa_2$ are constant. It is clear that
  $s=s(t)$ is in the form of the second formula in (\ref{Geody4}) if $s(0)=0$,
  and $P(s)$ is in the form of the second formula in
  (\ref{Geod003}) by plugging the above $s=s(t)$ into (\ref{Geody1}).
\end{ex}

\section{Construction of sprays from geodesics}

Given a family of curves $\mathcal{G}$ on a manifold, if
$\mathcal{G}$ can constitute the geodesics of a spray {\bf G} on
$M$, how can we solve {\bf G} (at least locally)? A spray induces
a (local) semispray and two pointwise projectively related sprays
induce a same semispray (cf. \cite{Shen6}). A semispray can also
be considered as a special parameterized family of curves, which
forms a path space.

In this section, we will start from a path space and introduce
some ways to construct sprays based on a path space and its
parameterization. We call it the geodesic method of construction
of sprays.

Similarly to a spray, a path space $\mathcal{G}$  on a manifold
$M$ is usually defined on a conical region $\mathcal{C}$ of
 $TM\setminus \{0\}$ (see Definition \ref{DGeod00}), and  $\mathcal{G}$ is called singular if
 $\mathcal{C}\ne
 TM\setminus \{0\}$.

\begin{Def}\label{DGeod00}
Let $\mathcal{G}$ be a family of $C^{\infty}$ parameterized curves
(called paths) on an $n$-dimensional manifold $M$.  $\mathcal{G}$
or $(M,\mathcal{G})$ is called an $n$-dimensional path space if on
a conical region $\mathcal{C}$ of $TM$ it satisfies

 {\rm (i)}  for $y\in \mathcal{C}_x$, there is a curve
 $\sigma:(-\epsilon,\epsilon)\rightarrow M$ in $\mathcal{G}$ with
  $\sigma'(0)=y$;

{\rm (ii)} for any  $\sigma,\tau$ in $\mathcal{G}$ with
$\sigma'(0)=\tau'(0)$,  $\sigma$ and $\tau$ coincide in  a small
intervals of $0$;

 {\rm (iii)} if a curve $\sigma$ is in $\mathcal{G}$, then for any constants $\lambda>0$ and $t_o$,
 the curve $\eta$ is also in
 $\mathcal{G}$, where $\eta$ is defined by $\eta(t):=\sigma(\lambda
 t+t_o)$.
\end{Def}

An equivalent version of Definition \ref{DGeod00} in regular case
is refereed to \cite{Shen6} ($P_{52}$).

\begin{ex}\label{EGeody1}
 Consider a set $\mathcal{G}$ of a family of curves $x=x(s)$ on $R^2$ in the form
 \beqn
 &&x(s)=\sigma(s;x_o,y_o),\ \ \ \ \big(x(0)=x_o=(a,b),\ \
 x'(0)=y_o=(u,v)\big),\\
 &&\sigma(s;x_o,y_o):=(a,b)+(u,v)s-(0,1)\big(\frac{1}{3}u^3s^3+au^2s^2\big),
 \eeqn
 where $a,b,u,v$ are arbitrary parameters. It can be directly verified that $\mathcal{G}$ is a path space on
 $R^2$, since  Definition
\ref{DGeod00} (i) (ii)  automatically hold, and Definition
\ref{DGeod00} (iii)  follows from
 $$
 \sigma(\lambda s+s_o;x_o,y_o)=\sigma(s;\hat{x}_o,\hat{y}_o),
 $$
 where we define
  \beqn
&&x_o=(a,b),\ \ \ \ y_o=(u,v),\ \ \ \
\hat{x}_o=(\hat{a},\hat{b}),\ \ \ \
\hat{y}_o=(\hat{u},\hat{v}),\\
&&\hat{a}:=a+us_o,\hspace{1.3cm}
\hat{b}:=b+vs_o-\frac{1}{3}u^2(3a+us_o)(s_o)^2,\\
&&\hat{u}:=\lambda u,\hspace{2cm} \hat{v}:=\lambda v-\lambda
u^2s_o(2a+us_o).
  \eeqn
\end{ex}

For a path space $\mathcal{G}$, we have different ways to
parameterize the paths in $\mathcal{G}$ under a parametric
variable and some constant parameters (see Theorem \ref{LGeody3}
and Lemma \ref{LGeody1}). Example \ref{EGeody1} satisfies
(\ref{Geodj0}) and (\ref{Geodj1}) in the following Lemma
\ref{LGeody1} with
 $$
f(s;x_o,y_o):=-(0,1)\big(\frac{1}{3}u^3s^3+au^2s^2\big),\ \
x_o=(a,b),\ y_o=(u,v).
 $$

 \begin{lem}\label{LGeody1}
 An $n$-dimensional path space $(M,\mathcal{G})$ is locally
 expressed as the following family of curves $x=x(s)$ with arbitrary constant parameters $x_o,y_o\in
 R^n$:
 \be\label{Geodj0}
 x(s)=\sigma(s;x_o,y_o)=x_o+y_os+f(s;x_o,y_o),
 \ee
where $f$ is a smooth function satisfying $f(0;x_o,y_o)=
f'(0;x_o,y_o)=0$ and
 \beq
&&f(s;\hat{x}_o,\hat{y}_o)=f(\lambda
s+s_o;x_o,y_o)-f(s_o;x_o,y_o)-\lambda f'(s_o;x_o,y_o)s,\label{Geodj1}\\
&&\ \ \ \big(\hat{x}_o:=x_o+y_os_o+f(s_o;x_o,y_o),\ \
\hat{y}_o:=\lambda y_o+\lambda f'(s_o;x_o,y_o)\big).\nonumber
 \eeq
\end{lem}

 It is clear that the collection of geodesics of a spray
naturally forms a path space.  Shen proves the converse in the
following lemma (\cite{Shen6}). We also give the proof for
convenience.

\begin{lem}\label{LGeody2}
 A path space $\mathcal{G}$ induces a spray {\bf G} with the
 set of geodesics of {\bf G} being $\mathcal{G}$.
\end{lem}

{\it Proof :} Let $(\mathcal{G},M)$ be a path space on a conical
region $\mathcal{C}$. For a given $y\in \mathcal{C}_x$, there is a
curve $\sigma:(-\epsilon,\epsilon)\rightarrow M$ in $\mathcal{G}$
with
  $\sigma(0)=x,\sigma'(0)=y$ by Definition
\ref{DGeod00} (i). Define
   $$
 G^i(y):=-\frac{1}{2}\frac{d^2\sigma}{ds^2}(0),
   $$
which is independent of the choice of $\sigma$ by Definition
\ref{DGeod00} (ii). We are going to verify that {\bf G} is a
spray. For any constant $\lambda>0$, let $\eta(s):=\sigma(\lambda
s)\in\mathcal{G}$ (see Definition \ref{DGeod00} (iii)). Then we
have
 $$
 G^i(\lambda
 y)=-\frac{1}{2}\frac{d^2\eta}{ds^2}(0)=-\frac{1}{2}\lambda^2\frac{d^2\sigma}{ds^2}(0)=\lambda^2G^i(y),
 $$
which implies that $G^i$ is positively homogeneous of degree two.
Further, for any $\eta:(a,b)\rightarrow M$ in $\mathcal{G}$ and
any fixed $t\in (a,b)$, define $\gamma(s):=\eta(s+t)$. Then  we
have
 $$
 \eta'(t)=\gamma'(0),\ \ \ \eta''(t)=\gamma''(0).
 $$
So by the definition of $G^i$, we get
 $$
 G^i\big(\eta'(t)\big)=G^i\big(\gamma'(0)\big)=-\frac{1}{2}\frac{d^2\gamma^i}{ds^2}(0)
 =-\frac{1}{2}\frac{d^2\eta^i}{ds^2}(t),
 $$
which implies that $\eta$ satisfies the following ODE:
 $$
 \frac{d^2\eta^i}{ds^2}+2G^i\big(\frac{d\eta}{ds}\big)=0.
 $$
Therefore, {\bf G} is a spray, and the set of geodesics of {\bf G}
coincides with $\mathcal{G}$.   \qed

\

In Lemma \ref{LGeody2}, by different choices of the parametric
variables, it (locally) induces a projective class $Proj({\bf G})$
of {\bf G}, each of which is projective to {\bf G}.

For a given path space $\mathcal{G}$, it induces a spray {\bf G}
by Lemma \ref{LGeody2}. Then {\bf G} defines a semispray
$\hat{\mathcal{G}}$ (see \cite{Shen6}: $P_{37}$) and the geodesics
of {\bf G} and $\hat{\mathcal{G}}$ are closely related (see
\cite{Shen6}: Lemma 3.1.1). Therefore, in $\mathcal{G}$, any path
can be locally expressed as
 $$
 x^a=x^a(x^1;u,v),\ \ \ (u,v\in R^{n-1},\ 2\le a\le n),
 $$
 where $u,v$ are free constant parameters. So all paths in an $n$-dimensional path space
 depend only on $2(n-1)$ free constant parameters, where the Jaccobi determinant is not zero, namely,
  $$
  \det
 \begin{pmatrix}
  \pa x^a/\pa u  & \pa x^a/\pa v \\
\pa y^a/\pa u  & \pa y^a/\pa v
\end{pmatrix}
\ne 0,\ \ \ \ \big(y^a:=\frac{dx^a}{dx^1}\big).
$$
Then  we obtain Theorem \ref{LGeody3} for the construction of
sprays based on the parametric equations of path spaces.

 If we write (\ref{Geodj2}) in the form
  \be\label{Geodj02}
 x(t)=\sigma(\lambda t+\mu,u,v),
  \ee
  where $\lambda,\mu$ are constant numbers, then under the $2n$
  constant parameters $\lambda,\mu,u,v$, this family of curves
  satisfies Definition  \ref{DGeod00} (i)(ii)(iii). For instance,
  the 2-dimensional path space in Example \ref{EGeody1} can be written as the
  following family of curves
   \beqn
 &&x(s)=\tau(\lambda s+\mu;b_o,v_o)=(0,b_o)+(1,v_o)(\lambda
 s+\mu)-\frac{1}{3}(0,1)(\lambda s+\mu)^3,\\
 &&\hspace{1cm}\big(\lambda:=u,\ \mu:=a,\ v_o:=\frac{v}{u}+a^2,\
 b_o=b-av_o+\frac{1}{3}a^3\big).
   \eeqn

   By Theorem \ref{LGeody3}, if the set $\mathcal{A}$ of a family of curves on an
   $n$-dimensional manifold defines a path
   space, then $\mathcal{A}$ just depends on $2(n-1)$ free
   constant parameters. For example, in $R^n$, all circles with
   fixed radius cannot define a path space when $n\ge 3$, because
   in this case, the circles depend on more than $2(n-1)$ free
   constant parameters.

\

Now we introduce a method of constructing a spray {\bf G}
determined by a path space considered as the geodesics of {\bf G},
which is similar to Okubo's method for the construction of a
Finlser metric from a hypersurface as its indicatrix. We can start
from a family of curves given by (\ref{Geodj0}) or (\ref{Geodj2})
to determine a corresponding spray.

\

\noindent {\bf Method (I)}: For a family of curves  given by
(\ref{Geodj0}) satisfying (\ref{Geodj1}), actually we can reduce
one constant parameter since (\ref{Geodj0}) can be written as (if
$y^1_o\ne 0$)
 \beqn
 &&x(t)=\sigma(t;x_o,\bar{y}_o)=x_o+\bar{y}_ot+f(t;x_o,\bar{y}_o),\\
 &&\ \ \ \ \big(t:=y^1_os,\ \bar{y}_o^a:=y^a_o/y^1_o,\ \bar{y}_o:=(\bar{y}_o^a)\big).
 \eeqn
Let a path space be determined by (\ref{Geodj0}) and we put
 \beq
&&x=x_o+y_os+f(s;x_o,y_o),\ \ \ y(=\frac{dx}{ds})=y_o+f'(s;x_o,y_o)),\label{Geodj3}\\
&&G^i:=-\frac{1}{2}\frac{d^2x}{ds^2}=-\frac{1}{2}f''(s;x_o,y_o)).\label{Geodj4}
 \eeq
 Then  we obtain a spray {\bf G} from (\ref{Geodj4})  by eliminating $x_o,y_o,s$ in (\ref{Geodj4}) from
 (\ref{Geodj3}), where $s$ is a geodesic parameter of the spray
 {\bf G}.

 \

 \noindent {\bf Method (II)}:  Suppose that a family of
 curves are given by the parametric equation (\ref{Geodj2}) with $2(n-1)$ free constant parameters
 $u,v$. This case is more convenient to construct sprays.
 With an auxiliary parameter $c>0$, we put
 \be\label{Geod0001}
 x=\sigma(cs;u,v),\ \ \ \
 y=\frac{dx}{ds}=c\frac{d\sigma}{d\hat{s}}(cs;u,v),\ \ \ \hat{s}:=cs.
 \ee
Theoretically, we can express $c,s,u,v$ as functions of $x,y$ from
(\ref{Geod0001}). Then plugging them into the following
 \be\label{Geod0002}
 G^i:=-\frac{1}{2}\frac{d^2x^i}{ds^2}=c^2\frac{d^2\sigma^i}{d\hat{s}^2}(cs;u,v),
 \ee
we obtain a spray {\bf G} given by (\ref{Geod0002}), where $s$ is
a geodesic parameter of the spray
 {\bf G}.

\

Now in the following Examples \ref{EGeody001}-\ref{EGeody004}, we
use Method (I) or Method (II) to show how we construct sprays from
given path spaces by eliminating the corresponding parameters.

\begin{ex}\label{EGeody001}
 Consider a set $\mathcal{G}$ of a family of curves on $R^3$:
  \beqn
 &&x(s)=(a,b,c)+(u,v,w)s-(0,1,0)h(s),\\
&& \big(h(s):=-\frac{1}{3}(u^3+w^3)s^3-(au^2+cw^2)s^2\big),
  \eeqn
  where $a,b,c,u,v,w$ are constant parameters. $\mathcal{G}$ is a
  path space. By (\ref{Geodj3}) we get
   \be\label{Geodj5}
    x^1=a+us,\ \ x^3=c+ws,\ \   y^1=u,\ \ y^3=w.
   \ee
   By (\ref{Geodj4}), the induced spray {\bf G} is given by
   \beqn
 G^1&&\hspace{-0.6cm}=-\frac{1}{2}\frac{d^2x^1}{ds^2}=0,\ \ \  \ \
 \ \
 G^3=-\frac{1}{2}\frac{d^2x^3}{ds^2}=0,\\
  G^2&&\hspace{-0.6cm}=-\frac{1}{2}\frac{d^2x^2}{ds^2}=(u^3+w^3)s+(au^2+cw^2)\\
  &&\hspace{-0.6cm}=(u^3+w^3)s+\big[(x^1-us)u^2+(x^3-ws)w^2\big] \ \ \big(by \ (\ref{Geodj5})\big)\\
&&\hspace{-0.6cm}=x^1u^2+x^3w^2=x^1(y^1)^2+x^3(y^3)^2\ \ \big(by \
(\ref{Geodj5})\big).
   \eeqn
    {\bf G} has zero Riemann curvature and so it is metrizable (a Finsler spray) (\cite{Yang2}).
\end{ex}

\begin{ex}\label{EGeody002}
 Let $\mathcal{G}$ be  the set of all circles with fixed radius $r$ on
 $R^2$. We parameterize $\mathcal{G}$ by
  $$
 x^1(s)=a+r\cos s,\ \ \ \ x^2(s)=b+r\sin s,
  $$
  where $a,b$ are arbitrary constant parameters. $\mathcal{G}$ depends on just two free constant parameters. By Theorem \ref{LGeody3},
  $\mathcal{G}$ defines a spray {\bf G} on $R^2$ with $s$ as a geodesic parameter of {\bf G}.
  We show the spray as follows. With an auxiliary parameter $c>0$,
  it follows from (\ref{Geod0001}) that
   \be\label{Geodj6}
 x^1=a+r\cos cs,\ \  x^2=b+r\sin cs,\ \  y^1=-cr\sin cs,\ \
 y^2=cr\cos cs.
 \ee
 Then plugging the latter two formulas of (\ref{Geodj6})  into (\ref{Geod0002})
 yields a spray {\bf G} given by
  \beqn
 &&G^1=-c^2r\cos cs=-\frac{1}{r}y^2\sqrt{(y^1)^2+(y^2)^2},\\
 &&G^2=-c^2r\sin
 cs=\frac{1}{r}y^1\sqrt{(y^1)^2+(y^2)^2}.
   \eeqn
   This circle spray first appears in \cite{Shen6} ($P_{49}$), and
 even locally it is not metrizable (\cite{Yang2}).
\end{ex}

\begin{ex}\label{EGeody003}
 Consider a family of semicircles $\mathcal{G}$ on the positive semi-plane $R^2_+$
 with center on $x^1$-axis and arbitrary radius. Note that $\mathcal{G}$ is singular
 at the direction parallel to $x^2$-axis. We can
 parameterize $\mathcal{G}$ by
  $$
 x^1=a+b\cos s,\ \ \ x^2=b\sin s,\ \ \ (x^2>0,\ b\ge 0),
  $$
   where $a,b$ are arbitrary constant parameters. $\mathcal{G}$ depends on just two free constant parameters.
   By Theorem \ref{LGeody3},
  $\mathcal{G}$ defines a spray {\bf G} on $R^2_+$ with $s$ as a geodesic parameter of {\bf
  G}. With an auxiliary parameter $c>0$,
  by (\ref{Geod0001}) we get
   \be\label{Geodj06}
 x^1=a+b\cos cs,\ \ \ x^2=b\sin cs,\ \ \ y^1=-bc\sin cs,\ \ \
 y^2=bc\cos cs.
   \ee
   Then similarly, by the elimination of the parameters $a,b,c,s$ in (\ref{Geod0002}) from
   (\ref{Geodj06}),
     the spray {\bf G} with $s$ being a geodesic parameter is given
   by
    \be\label{Geodj0006}
 G^1=-\frac{y^1y^2}{2x^2},\ \ \ \  G^2=\frac{(y^1)^2}{2x^2}.
    \ee
The spray {\bf G} is regular on $R^2_+$ (any straight lines
parallel to $x^2$-axis are geodesics of {\bf G}). {\bf G} is of
isotropic curvature, and locally it is not metrizable by the
method in \cite{BM2, Yang2}.
\end{ex}

\begin{ex}\label{EGeody004}
 Let $B^n$ be the unit ball in $R^n$ and $\mathcal{G}$ be all
 circle arcs in $B^n$ which are perpendicular to the boundary
 $S^{n-1}=\pa B^n$. Let $s$ be the arc-length parameter of a
 circle arc induced by the Euclidean metric.  What is the spray {\bf G} induced by $\mathcal{G}$
 with $s$ being a geodesic parameter of {\bf G} (see Example 4.1.4 in \cite{Shen6})? We will show that  {\bf G} is given by
  \be\label{Geodj7}
 G^i=\frac{\langle x,y\rangle y^i-|y|^2x^i}{1-|x|^2},
  \ee
which is not metrizable by \cite{Yang2}.
 Now for arbitrarily given $p,q\in S^{n-1}$, there is a circle arc
  $\gamma$ in $\mathcal{G}$, in which $\gamma$ is perpendicular to
  $S^{n-1}$ at $p,q$. Let $C$ be the circle with $\gamma \subset
  C$. The center and radius of $C$ are respectively given by
   $$
 \tau(p+q),\ \ \ \ |p-\tau(p+q)|,\ \ \  (\tau:=(1+pq)^{-1}),
   $$
   where $pq$ is the Euclidean inner product of $p,q$. Then
   $\gamma$ is parameterized by the equation
    \be\label{Geodj8}
 x(s)=x(s;p,q)=[p-\tau(p+q)]\cos s+|p-\tau(p+q)|p\sin s+\tau(p+q).
    \ee
    Since there are just $2(n-1)$ free constant parameters in
    (\ref{Geodj8}), the family of curves in the form
    (\ref{Geodj8}) define a path space by Theorem \ref{LGeody3}. Now based on (\ref{Geod0001}) and
    (\ref{Geod0002}), we can give the spray {\bf G} from (\ref{Geodj8})
    with $s$ being a geodesic parameter of {\bf G}.  With an auxiliary parameter $c>0$, by
    (\ref{Geod0001}) we put
     \beq
  &&x=[p-\tau(p+q)]\cos cs+|p-\tau(p+q)|p\sin cs+\tau(p+q),\label{Geodj9}\\
  &&y=-c[p-\tau(p+q)]\sin cs+c|p-\tau(p+q)|p\cos cs.\label{Geodj10}
     \eeq
  By (\ref{Geod0002}) we have
   \be\label{Geodj11}
 2G^i:=c^2\big\{[p-\tau(p+q)]^i\cos cs+|p-\tau(p+q)|p^i\sin
 cs\big\}.
   \ee
   By a direct lengthy computation, we can eliminate the parameters
   $p,q,c,s$ in (\ref{Geodj11}) from (\ref{Geodj9}) and
   (\ref{Geodj10}) (the details are omitted). Finally, the spray
   {\bf G} is given by (\ref{Geodj7}).
\end{ex}

\section{Projectively complete sprays}

For a given spray {\bf G}, if we know the general solutions of all
geodesics of {\bf G}, then under another parameter as a geodesic
parameter, we can determine a corresponding spray projectively
related to {\bf G}. Now suppose that the general solutions of
geodesics of {\bf G} are locally given by
 \be\label{Geod00}
 x=\sigma(t)=\sigma(t;u,v),\ \ \ \ \big(u,v\in R^{n-1}\big),
 \ee
where $t$ is a geodesic parameter of {\bf G} and $u,v$ are free
constant parameters. Sometimes, it is also convenient to put
$u=\sigma(0),v=\sigma'(0)$  for the elimination of parameters.
Make a change of the variables from $t$ to $s$ with
 \be\label{Geod000}
 t=t(s)=t(s;u,v),\ \ \ \frac{dt}{ds}>0.
 \ee
With an auxiliary parameter $c>0$, we put
 \be\label{Geod01}
 x=\sigma(t(cs);u,v),\ \ \ \
 y=\frac{dx}{ds}=c\frac{d\sigma}{dt}\frac{dt}{ds},
 \ee
 where $dt/ds$, as a function of $s$, takes the value at $cs$. Further, we have
 \beq
 \frac{d^2x^i}{ds^2}&&\hspace{-0.6cm}=c^2\frac{d^2\sigma^i}{dt^2}\big(\frac{dt}{ds}\big)^2
 +c^2\frac{d\sigma^i}{dt}\frac{d^2t}{ds^2}\nonumber\\
&&\hspace{-0.6cm}=-2G^i(x,\frac{d\sigma}{dt})c^2\big(\frac{dt}{ds}\big)^2
+c^2\frac{d\sigma^i}{dt}\frac{d}{dt}\big(\frac{dt}{ds}\big)\frac{dt}{ds}\nonumber\\
&&\hspace{-0.6cm}=-2G^i(x,y)+c\frac{d}{dt}\big(\frac{dt}{ds}\big)y^i.\label{Geod02}
 \eeq
Expressing $c,t$ in terms of $x,y$ from (\ref{Geod01}), and then
plugging $c,t$ into (\ref{Geod02}), we obtain a spray ${\bf
\bar{G}}$ given by
 \beq\label{Geod03}
 \bar{G}^i&&\hspace{-0.6cm}=G^i-\frac{1}{2}\frac{d}{dt}\big(\frac{dt}{ds}\big)cy^i=G^i+Py^i,\\
&&\hspace{-0.6cm}\Big(P=P(x,y):=-\frac{1}{2}\frac{d}{dt}\big(\frac{dt}{ds}\big)c\Big),\nonumber
 \eeq
with $s$ being a geodesic parameter of ${\bf \bar{G}}$.

\begin{lem}
Suppose that the general solutions of geodesics of a spray {\bf G}
are given by (\ref{Geod00}). Let $s$ be another parameter related
to $t$ by (\ref{Geod000}). Then a spray ${\bf \bar{G}}$ projective
to {\bf G} with $s$ being its geodesic parameter is given by
(\ref{Geod03}), where $c,t$ are determined by (\ref{Geod01}).
\end{lem}

Under certain condition, a spray can be projectively
(positively/negatively) complete, which is shown in Theorem
\ref{th3}. Now we give the proof of Theorem \ref{th3}.

\

{\it Proof of Theorem \ref{th3} :} Let {\bf G} be a spray on a
manifold $M$. For an arbitrary geodesic $x=x(t)$, suppose that $t$
belongs to the
 maximal interval $I$  given by (\ref{Geodc1}).

If $I=(a,+\infty)$ or
 $I=(-\infty,b)$, we respectively make a change of the variables
 from $t$ to $s$ by
  \be\label{Geodc2}
 s=\ln(1-\frac{t}{a}),\ \ \ or \ \ \ s=-\ln(1-\frac{t}{b}),
  \ee
  either of which gives $s(0)=0,s'(t)>0$ and the maximal interval of $s$
  with $s\in(-\infty,+\infty)$.

  If $I=(a,b)$, make a change by
  (\ref{Geodc2}) and then we respectively have
   $$
 s\in\big(-\infty,\ln(1-\frac{b}{a})\big),\ \ \ or \ \ \
 s\in\big(-\ln(1-\frac{a}{b}),+\infty).
   $$

   If $I=(a,b)$, make a change of the variables
 from $t$ to $s$ by
 \be\label{Geodc3}
 s=\ln\frac{1-t/a}{1-t/b},\ \ \ or \ \ \ s=\tan\big[\frac{\pi}{b-a}(t-\frac{a+b}{2})\big]+\tan\big(\frac{b+a}{b-a}\frac{\pi}{2}\big),
  \ee
  either of which gives $s(0)=0,s'(t)>0$ and the maximal interval of $s$
  with $s\in(-\infty,+\infty)$.

  Therefore, by the change
  (\ref{Geodc2}) or (\ref{Geodc3}), we obtain a (positively/negatively) complete spray which is projective
  to {\bf G}. This completes the proof.   \qed

\

As an application of Theorem \ref{th3}, we give the following
Examples \ref{EGeod1}-\ref{EGeod3} to show the construction of the
(positively/negatively) complete sprays projective to given
sprays.

\begin{ex}\label{EGeod1}
Let $F$ be the Funk metric on a strongly convex domain
$\Omega\subset R^n$. The Minkowski spray ${\bf G}=0$ on $\Omega$
has its geodesics given by
 $$
 x(t)=vt+u,\ \ \ \big(-\frac{1}{F(u,-v)}<t<\frac{1}{F(u,v)}\big),
 $$
 where $u,v\in R^n$ are arbitrary constant
 vectors. By (\ref{Geodc2}), put $t=t(s)$ as
 \be\label{Geod0003}
s=-\ln\big[1-tF(u,v)\big].
  \ee
With   $s$ being a geodesic parameter, we obtain a projectively
flat and positively complete spray ${\bf \bar{G}}$, which  will be
shown to be the Finsler spray induced by
 $F$, namely,
  \be\label{Geod3}
 \bar{G}^i=\frac{1}{2}Fy^i.
  \ee

Actually, it follows from (\ref{Geod003}) and (\ref{cwy3}) that
   \be\label{Geod4}
 \frac{dt}{ds}=\frac{1}{F(u,v)}-t=\frac{1}{F(vt+u,v)}.
   \ee
Then (\ref{Geod01}) gives
   \be\label{Geod5}
 x=vt+u,\ \ \  y=cv\frac{dt}{ds}.
   \ee
    It is clear from (\ref{Geod5}) and (\ref{Geod4}) that
    \be\label{Geod7}
 F(x,y)=cF(vt+u,v)\frac{dt}{ds}=c.
    \ee
    Therefore, by (\ref{Geod03}), (\ref{Geod4}) and (\ref{Geod7}), the spray {\bf G} is given by
    (\ref{Geod3}).
\end{ex}

\begin{ex}\label{EGeod2}
 In Example \ref{EGeod1}, by (\ref{Geodc3}), put $t=t(s)$ as
  \be\label{Geod007}
 s=\ln\frac{1+tF(u,-v)}{1-tF(u,v)}.
  \ee
  With $s$ being a geodesic parameter, we obtain a projectively flat and complete spray
  ${\bf \bar{G}}$, which  will be
shown to be the Finsler spray induced by the Klein metric
 $\bar{F}(x,y):=\big[F(x,y)+F(x,-y)\big]/2$, namely,
   \be\label{Geod8}
 \bar{G}^i(x,y)=\frac{1}{2}\big[F(x,y)-F(x,-y)\big]y^i.
   \ee

   Firstly, by (\ref{Geod007}), we get
    \beq
 \frac{dt}{ds}&&\hspace{-0.6cm}
 =\frac{\big[1-tF(u,v)\big]\big[1+tF(u,-v)\big]}{F(u,v)+F(u,-v)},\label{Geod9}\\
  \frac{d}{dt}\big(\frac{dt}{ds}\big)&&\hspace{-0.6cm}=\frac{F(u,-v)-F(u,v)-2tF(u,v)F(u,-v)}{F(u,v)+F(u,-v)}.\label{Geod10}
    \eeq
    Secondly,  (\ref{Geod01}) gives
  $$
 x=vt+u,\ \ \
 y=cv\frac{dt}{ds},
   $$
from which we have
 \be\label{Geod12}
 F(x,y)=F(vt+u,v)c\frac{dt}{ds},\ \ \ \
  F(x,-y)=F(vt+u,-v)c\frac{dt}{ds}.
 \ee
Plugging (\ref{cwy3}), (\ref{cwy003}) and (\ref{Geod9}) into
(\ref{Geod12}), we obtain
 \be\label{Geod13}
 c=F(x,y)+F(x,-y),\ \ \
 t=\frac{F(x,y)F(u,-v)-F(x,-y)F(u,v)}{F(u,v)F(u,-v)[F(x,y)+F(x,-y)]}.
 \ee
Finally, by (\ref{Geod13}) and (\ref{Geod10}), it follows from
(\ref{Geod03}) that the spray ${\bf \bar{G}}$ is given by
(\ref{Geod8}).
\end{ex}

\begin{ex}\label{EGeod3}
 For the spray {\bf G} in Example \ref{EGeod1}, we will introduce a
 different way from that in Example \ref{EGeod2} to make {\bf G}
 be complete, which is actually to use (\ref{Geodc2}) to make
 complete the Finsler spray  induced by the Funk
 metric $F$. For a geodesic $x(t)=vt+u$ of {\bf G}, put
  $$
 s=\ln\Big[1-\frac{\ln(1-tF(u,v))}{a}\Big],\ \ \
 a:=\ln\Big[1+\frac{F(u,v)}{F(u,-v)}\Big].
  $$
  In a similar way to that for the computation in Example \ref{EGeod2}, we  obtain a projectively flat and complete spray
  ${\bf \bar{G}}$ given as follows:
  $$
  \bar{G}^i(y)=G^i_F(y)+\frac{1}{2}\frac{F(y)}{\ln\frac{F(-y)}{F(y)+F(-y)}}y^i,\
  \ \ \Big(G^i_F(y):=\frac{1}{2}F(y)y^i\Big).
  $$
${\bf \bar{G}}$  is of scalar curvature and actually we can verify
that ${\bf \bar{G}}$ is not metrizable by using  the method in
\cite{BM2}.
\end{ex}

\begin{ex}\label{EGeod4}
 For the family of semicircles $\mathcal{G}$ on $R^2_+$ as shown in Example \ref{EGeody003}, we can
 parameterize them in the following form
 \be\label{ygj1}
 x^1=u-v\sin t,\ \ \ x^2=v\cos t,\ \ \ (x^2>0,\ v\ge 0),
  \ee
   where $u,v$ are arbitrary constant parameters. We have shown in
Example \ref{EGeody003} that the spray {\bf G} determined by
$\mathcal{G}$ is given by (\ref{Geodj0006}), that is,
 \be\label{ygj2}
 G^1=-\frac{y^1y^2}{2x^2},\ \ \ \  G^2=\frac{(y^1)^2}{2x^2}.
 \ee
 We can make {\bf G} be projectively complete on the conical
 region $\mathcal{C}$ with the direction $(0,1)$ being deleted from $TR^2_+\setminus
 \{0\}$. Since $-\pi/2< t<\pi/2$ for any $u,v$ in (\ref{ygj1}), by
 (\ref{Geodc3}), we let
  $$
 s=\tan t.
  $$
  Then by (\ref{Geod03}), we get a complete spray ${\bf \bar{G}}$
  projective to {\bf G} with the projective factor $P$ being given by
   \be\label{ygj3}
  P=c\Big[-\frac{1}{2}\frac{d}{dt}\big(\frac{dt}{ds}\big)\Big]_{t=cs}=c\big[\cos t\sin
  t\big]_{t=cs}=\frac{c^2s}{1+c^2s^2}.
   \ee
   Now it follows from (\ref{Geod01}) that
    \beqn
 && x^1=u-\frac{vcs}{\sqrt{1+c^2s^2}},\hspace{1cm}
  x^2=\frac{v}{\sqrt{1+c^2s^2}},\\
  &&y^1=\frac{-vc}{(1+c^2s^2)^{3/2}},\hspace{1.25cm}
  y^2=\frac{-bc^2s}{(1+c^2s^2)^{3/2}},
    \eeqn
    from which we get
    $$
  s=-\frac{x^2y^2}{(y^1)^2+(y^2)^2},\hspace{1cm}
  c=-\frac{(y^1)^2+(y^2)^2}{x^2y^1}.
    $$
    Plugging $s,c$ in the above into (\ref{ygj3}) yields
 $P=-y^2/x^2$. Thus the spray ${\bf \bar{G}}$ is given by
  $$
 \bar{G}^1=G^1+Py^1=-\frac{3y^1y^2}{2x^2},\ \ \ \  \bar{G}^2=G^2+Py^2=\frac{(y^1)^2-2(y^2)^2}{2x^2}.
  $$
${\bf \bar{G}}$ is  complete on the conical region $\mathcal{C}$
but not complete in the direction $(0,1)$. We don't know whether
the spray {\bf G} in  (\ref{ygj2}) can be projectively complete or
not on $TR^2\setminus\{0\}$.  Besides, ${\bf \bar{G}}$ is of
isotropic curvature, and locally it is not metrizable by the
method in \cite{BM2, Yang2}.
\end{ex}

\vspace{0.5cm}

\noindent Guojun Yang \\
Department of Mathematics \\
Sichuan University \\
Chengdu 610064, P. R. China \\
yangguojun@scu.edu.cn

\end{document}